\documentclass[10pt]{extarticle}

\usepackage{amsmath,amssymb,amsthm}

\def\RR{{\mathbb R}}
\def\NN{{\mathbb N}}

\def\EE{{\mathbb E}}
    
\newcommand{\ba}{\begin{array}} \newcommand{\ea}{\end{array}}
\newfont{\bsl}{cmbxsl10 scaled 1095}

\newcommand{\eps}{\varepsilon}

\newfont{\deu}{eufm10 scaled 1000}

\newcommand{\co}{\operatorname{co}}

\newcommand{\beq}{\begin{equation}}
\newcommand{\eeq}{\end{equation}}

\begin{document}

\title{R.E. Bruck, proof mining and a rate of asymptotic regularity for 
ergodic averages in Banach spaces}

\author{Anton Freund and Ulrich Kohlenbach\\ 
Department of Mathematics\\ Technical University of Darmstadt \\
 Schloss\-garten\-str.~7\\ 64289~Darmstadt, Germany \\ 
\{freund,kohlenbach\}@mathematik.tu-darmstadt.de}

\date{June 30, 2022}

\maketitle 
\begin{abstract}
\noindent We analyze a proof of Bruck to obtain an explicit rate of asymptotic regularity for Ces\`aro means in uniformly convex Banach spaces. Our rate will only depend on a norm bound and a modulus~$\eta$ of uniform convexity. One ingredient for the proof by Bruck is a result of Pisier, which shows that every uniformly convex (in fact every uniformly nonsquare) Banach space has some Rademacher type~$q>1$ with a suitable constant~$C_q$. We explicitly determine~$q$ and~$C_q$, which only depend on the single value~$\eta(1)$ of our modulus. Beyond these specific results, we summarize how work of Bruck has inspired developments in the proof mining program, which applies tools from logic to obtain results in various areas of mathematics.
\end{abstract}

\section{Introduction} 
\newtheorem{definition}{Definition}
\newtheorem{proposition}[definition]{Proposition}
\newtheorem{remark}[definition]{Remark}
\newtheorem{theorem}[definition]{Theorem}
\newtheorem{corollary}[definition]{Corollary}
\newtheorem{lemma}[definition]{Lemma}
\newtheorem{exercise}[definition]{Exercise}
\newtheorem{clm}[definition]{Claim}
\newtheorem{prop}[definition]{Proposition}
\newtheorem{example}[definition]{Example}
\newtheorem{notation}[definition]{Notation}
\newtheorem{application}[definition]{Application}

\maketitle

Proof mining is the  project of applying 
proof-theoretic transformations to obtain new quantitative and qualitative 
information from given proofs in areas of core mathematics 
such as nonlinear analysis, convex optimization and geodesic geometry 
(see e.g. \cite{Kohlenbach(ICM)}). Bruck, who himself did fundamental 
work on quantitative issues of metric fixed point theory 
(\cite{BaillonBruck}), was a major source of inspiration in this program 
both by providing in his research deep results which naturally 
asked for a more finitary quantitative treatment as well as by 
introducing fundamental new notions which were particularly suited 
for such a proof-theoretic enterprise. 
\\ In Section \ref{section-survey} 
we will give a short survey on the important role
which results of Bruck have had in the development of proof mining.
\\ The most recent proof-theoretic analysis of a work of Bruck has been 
carried out in \cite{FreundKohlenbach} which gives an explicit rate of 
metastability (in the sense of T. Tao) 
for a strong nonlinear ergodic theorem in uniformly convex 
Banach spaces due to \cite{KobMiy}  which in turn is crucially based on Bruck's 
seminal work on the nonlinear mean ergodic theorem \cite{Bruck79} 
and the convex approximation property \cite{Bruck(81)}. In \cite{Bruck(81)},
the asymptotic regularity of Ces\`aro means in uniformly convex Banach 
spaces is established using that uniformly convex Banach spaces $X$, which are thus B-convex, 
have a nontrivial Rademacher type~$q>1.$ This latter fact 
was first established by Pisier in \cite{pisier73}. \\
In Section \ref{section-Pisier} we will extract from Pisier's proof 
explicit lower estimates $>1$ for $q$ and an upper estimate for the 
relevant constant $C_q$ witnessing 
that the space has Rademacher type $q$ in terms of a constant $\delta>0$ 
which witnesses that $X$ is uniformly nonsquare. If $X$ is uniformly 
convex and $\eta$ is some modulus of uniform convexity, then $\delta$ 
can be taken as $\delta:=\eta(1).$
\\
In Section \ref{section-Bruck} we will then give an explicit rate 
of asymptotic regularity in terms of $\eta, q,C_q$ and so - by Section 
\ref{section-Pisier} - in terms of $\eta$ alone. Note that we do not use here 
{\bf the} optimal modulus $\delta_X$ of uniform convexity but any 
function $\eta$ witnessing the $\forall\eps\exists \delta$-definition 
of uniform convexity. In Remark~\ref{rmk:zhu} we will comment on an alternative
way to obtain a rate of asymptotic regularity, which combines our work with a result of
Zhu, Huang and Li~\cite{zhu}.

\section{Bruck and proof mining}\label{section-survey}
The work of Bruck contains many important results which he proved by 
prima facie noneffective means but which do ask - e.g. by their general 
logical form - for additional computational information. One such 
example is his convergence result 
on an iteration scheme (for suitable sequences $(\lambda_n),(\theta_n)$ 
in $[0,1]$)
\[ x_{n+1}=\left(1-\lambda_n\right)x_n+\lambda_nTx_n-\lambda_n\theta_n
\left(x_n-x_1\right) \]
for demicontinuous (single valued) 
pseudo-contractions $T:C\to C$ on nonempty 
closed and convex subsets $C$ of a Hilbert space $X$
(\cite[Theorem 4]{Bruck(74)}). While the Mann iteration of $T$ is  
not even asymptotically regular for pseudo-contractions (see 
\cite{Mutan}) and the Ishikawa iteration is asymptotically 
regular but strongly convergent only for compact $C$, 
Bruck's hybrid scheme converges strongly towards a fixed point of $T,$ 
if $T$ has a fixed point, which e.g.\ is the case when $C$ additionally 
is bounded (Bruck's scheme has been studied also in Banach spaces, see e.g. 
\cite{Reich(77),Reich(78),Reich(80)}). By a proof-theoretic 
analysis of Bruck's noneffective proof, 
K\"ornlein extracted in \cite{Koernlein(diss),Koernlein16} an explicit 
effective rate of metastability in the sense of T. Tao for this convergence 
result. This work is based in turn on \cite{KoernleinKohlenbach11}, where the special case 
of Lipschitzian pseudo-contractions had been treated. In this situation 
one has the asymptotic regularity of the iteration, i.e.
\[ \lim_{n\to\infty} \| x_n -Tx_n\| =0, \] 
even in general Banach spaces and with easier conditions on the scalars 
involved by a result from \cite{ChidumeZegeye} from which a polynomial 
rate of asymptotic regularity is extracted in \cite{KoernleinKohlenbach11}.
\\[1mm] 
Another inspiration for proof mining came from 
Bruck's nonconstructive proof of the existence of sunny
nonexpansive retractions \cite{Bruck(73)} onto fixed point sets $\operatorname{Fix}(T)$ 
of nonexpansive mappings $T:C\to C$ with $C\subseteq X$, where 
$X$ is reflexive, $C$ is bounded, closed and convex 
(actually Bruck's conditions are still weaker) and $T$ satisfies the so-called 
conditional fixed point property CFP. These conditions always hold e.g. for 
uniformly smooth spaces.
In \cite{Bruck(70)}, Bruck had shown that 
$\operatorname{Fix}(T)$ is a nonexpansive retract whenever the underlying space $X$ is 
a reflexive strictly convex Banach space making use of Zorn's lemma. 
If $X$ is uniformly smooth, then there is even a sunny nonexpansive 
retraction and there 
can -- already in smooth spaces -- only be at most one such sunny 
nonexpansive retraction~\cite{Bruck(73A)}. 
The first more constructive way of approaching this unique 
sunny nonexpansive retraction
is Reich's~\cite{Reich(80)} fundamental result that in 
uniformly smooth Banach spaces the sequence $(x_n)$ defined 
by $x_n$ being the fixed point of the strict contraction
\[ T_n:C \to C\quad\text{with}\quad T_n(x):= \left( 1-\frac{1}{n}\right) T(x)+\frac{1}{n}u\quad\mbox{for}\quad u\in C \]
strongly converges to the sunny nonexpansive 
retraction onto $\operatorname{Fix}(T)$ applied to $u.$ Reich's result actually is 
more general and, in particular, even applies to 
the aforementioned pseudo-contractions. However, the convergence proof 
is highly noneffective and, in fact, one can show that even for simple 
situations (like $X:=\RR$) there in general is no computable rate 
of convergence. 
\\ 
In \cite{KohlenbachSipos}, a completely constructive proof for the 
metastable reformulation of the convergence of $(x_n)$ (which -- 
noneffectively -- trivially is equivalent to the latter) is given 
together with an explicit rate of metastability.
\\[1mm] 
A different aspect by which Bruck's research is particularly attractive 
in the course of the proof mining paradigm is his formulation of 
important classes of 
mappings such as the firmly nonexpansive mappings 
\cite{Bruck(73)} and -- together with Reich -- the averaged and 
strongly nonexpansive (SNE) mappings \cite{BruckReich} as well as 
the strongly quasi-nonexpansive mappings \cite{Bruck(82)}. 
All these classes beautifully fit the requirements for
proof-theoretically well-behaved classes of nonlinear mappings 
and play a crucial rule in papers using the proof mining methodology such as
\cite{Ariza,Nicolae,Kohlenbach(SNE)}. 
For example, the condition of being firmly nonexpansive is purely universal and hence 
`tame' in the context of the functional interpretations used in proof mining.
The same is true for averaged mappings once one has an averaging 
constant $\alpha\in (0,1)$ together with a witness $N\in\NN$ such
that $\alpha\in [\frac{1}{N},1-\frac{1}{N}]$ given. In this case, 
extracted bounds will additionally depend on $N$ (but not on $\alpha$).
Being SNE in its original formulation seemingly has a very high 
logical complexity but, in fact, is equivalent to the existence 
of a number-theoretic function $\omega$ satisfying a purely universal 
condition (see \cite{Kohlenbach(SNE)}): 
\[ \ \ba{l} \forall c,k\in\NN \,\forall x,y\in C \\
\hspace*{1cm} \left( \| x-y\|\le c\wedge 
\| x-y\|-\| Tx-Ty\| <2^{-\omega(k)} \to \| (x-y)-(Tx-Ty)\| \le 2^{-k}
\right). \ea \] 
Then extractable bounds will depend additionally on $\omega.$ 
Moreover, from any proof of the property SNE for a class of mappings which 
satisfy the logical conditions in the metamathematical bound extraction 
theorems in proof mining one can extract such an SNE-modulus. This e.g.\ 
has been done for the class of firmly nonexpansive mappings in uniformly 
convex Banach spaces in \cite{Kohlenbach(SNE)} and for the averaged 
mappings in Hilbert space in \cite{Sip(20)}. The important property of 
SNE-mappings being closed under composition results in a simple computation 
of an SNE-modulus for compositions of SNE-mappings in terms 
of SNE-moduli for the individual mappings.
All this 
plays a crucial role in the extraction of a polynomial rate of 
asymptotic regularity in \cite{Kohlenbach(inconsistent)} for 
Bauschke's solution of the zero displacement conjecture (see also 
the recent generalization of \cite{Kohlenbach(inconsistent)} in 
\cite{Sip(20)}) as well as the quantitative analysis of proximal point 
type algorithms in \cite{Kohlenbach(PPA),Kohlenbach(HPPA),Kohlenbach(PPAcomonotone)}.
In the latter papers it is crucially used that all firmly nonexpansive 
mappings in uniformly convex Banach spaces (and so all resolvents) 
have a common SNE-modulus and that all averaged mappings in Hilbert 
space with some control on the averaging constants have a common 
SNE-modulus as well. \\ Also the property of being strongly quasi-nonexpansive
is logically very well-behaved (when localized to some fixed point $p$ of 
$T$) and gives rise to a corresponding modulus which 
witnesses this property quantitatively and which is used in the extraction 
of explicit rates of asymptotic regularity and metastability of 
algorithms which compute common fixed points of such mappings in 
geodesic settings such as CAT$(\kappa)$-spaces with $\kappa> 0$ 
(see \cite{Kohlenbach(SNE)}).
\\[1mm]
As mentioned already in the introduction, in the course of extracting 
explicit rates for a strong nonlinear ergodic theorem due to \cite{KobMiy}
which is carried out in \cite{FreundKohlenbach} 
we recently, in particular, analysed two proofs from \cite{Bruck79} 
and \cite{Bruck(81)}, respectively. The first proof concerns the existence of a 
convex continuous and unbounded function $\gamma:[0,\infty)\to [0,\infty)$ 
such that $\gamma(0)=0$ and for all $x_1,x_2\in C$ and $\lambda \in [0,1]$
\[ \gamma(\| T(\lambda x_1+(1-\lambda)x_2)-(\lambda Tx_1+(1-\lambda)Tx_2)\|) 
\le \| x_1-x_2\| -\| Tx_1 -Tx_2\| \]
and the second one proves the convex approximation property of B-convex 
Banach spaces using the fact 
that such spaces have Rademacher type $q>1$ by a result due to Pisier 
\cite{pisier73}. 
The latter property is used in \cite{Bruck(81)} to show that in uniformly 
convex Banach spaces the Ces\`aro means 
\[ x_n:=\frac{1}{n}\sum^{n-1}_{i=0} T^ix\] 
of a nonexpansive map~$T$ satisfy the asymptotic regularity property 
\[ \lim_{n\to\infty} \| x_n-Tx_n\| = 0.\]
In the remaining sections we will provide an explicit rate of convergence that only depends on some 
norm bound $b$ and a given modulus of uniform convexity $\eta$.
   
\section{On the Rademacher type of uniformly nonsquare Banach spaces}\label{section-Pisier}

In this section, we analyse a proof of Pisier~\cite{pisier73} to extract quantitative information on the Rademacher types of uniformly nonsquare Banach spaces. We then derive information on a probabilistic characterization of $B$-convexity, which is also due to Pisier~\cite{pisier73}. This supplements the quantitative analysis of a proof by Bruck~\cite{Bruck(81)} that has been carried out in~\cite{FreundKohlenbach}.

Consider a sequence $(\eps_i)$ of independent random variables that take values~$\pm 1$ with probability~$1/2$. These are conveniently realized by Rademacher functions~$r_i$ on $[0,1]$ (see e.\,g.~\cite[Theorem~2.b.3]{Lindenstrauss-Tzafriri-Part1}). For $1\le q<\infty$ and points $x_1,\ldots,x_n$ in a given Banach space~$X$, we consider the expected value
\[\EE\left(\left\| \sum^n_{i=1}\eps_ix_i\right\|^q\right)=\frac{1}{2^n}
\sum_{\eps_1,\ldots,\eps_n\in\{ -1,1\}}\left\| \sum^n_{i=1}\eps_ix_i\right\|^q
=\int^1_0\left\| \sum^n_{i=1} r_i(t)x_i\right\|^q dt.\]
We now recall the definition of Rademacher type (see~\cite[Section~9.2]{ledoux-talagrand} or \cite[Definition~1.e.12]{Lindenstrauss-Tzafriri}):
\begin{definition}\label{def:Rademacher-type}
A Banach space $(X,\|\cdot\|)$ is said to have Rademacher type 
$q\in [1,2]$ with constant~$C_q$ if all finite sequences $(x_1,\ldots,x_n)$ in $X$ validate
\[ \tag{$\ast$}\ \EE\left(\left\| \sum^{n}_{i=1} \eps_ix_i\right\|^q\right)^{1/q} 
\le C_q\cdot\left(\sum^n_{i=1} 
\| x_i\|^q\right)^{1/q}. \]
\end{definition}

Let us also recall the notion of uniformly nonsquare Banach space, which is due to James~\cite{James(64)} and can be equivalently expressed as follows:

\begin{definition}\label{nonsquare}
A Banach space $(X,\|\cdot\|)$ is uniformly nonsquare if there exists 
a $\delta \in (0,1]$ such that
\[ \min\left\{ \frac{\| x-y\|}{2},\frac{\| x+y\|}{2}\right\} \le (1-\delta)\cdot\max\{ \| x\|,\| y\|\}\]
holds for all $x,y\in X$.
\end{definition} 

We can now state our goal more precisely: Given~$\delta>0$, we want to determine~$q$ and $C_q$ such that any Banach space that is uniformly nonsquare for $\delta$ has Rademacher type $q$ with constant~$C_q$. Let us first give a quantitative version of the result that uniformly convex entails uniformly nonsquare:

\begin{lemma}\label{lem:convex-to-nonsquare}
Consider a uniformly convex Banach space $(X,\|\cdot\|)$ with a given modulus of convexity $\eta:(0,2]\to (0,1]$, which means that all $\eps\in(0,2]$ and $x,y\in X$ validate
\[ 
\| x\|,\|y\|\le 1\text{ and }\| x-y\|\ge \varepsilon\quad\Rightarrow\quad 
\left\| \frac{x+y}{2}\right\| \le 1-\eta(\varepsilon).\]
Then $X$ is uniformly nonsquare, where we may take 
$\delta:=\eta(1)$ in Definition~\ref{nonsquare}.
\end{lemma}
Also note that if Definition~\ref{nonsquare} holds for $\delta$, then $1-\delta$ is an upper bound for $\lambda_2(X)$ from \cite[p.~VII.1]{pisier73}.\\[1ex]
{\bf Proof} (see also the proof of 
Theorem 2.2.5 in \cite{Fackler(11)}):
Put $\lambda:=1-\eta(1)< 1.$ Since we may assume that $X$ is non-trivial,
let $e\in X$ be such that $\| e\| =1.$ Then the fact that $\eta$ is a 
modulus of convexity (applied to $x:=0$, $y:=e$ and $\eps:=1$) yields 
$\eta(1)\le \frac{1}{2}$ and thus $\lambda\in [\frac{1}{2},1)$. For any $x,y\in X$, we show
\[ \min\{\| x-y\|,\| x+y\|\}\le2\lambda\cdot\max\{ \| x\|,\| y\|\}.\]
Without loss of generality, we assume 
$\| x\|\ge \| y\|$. If we have $\| x\|=0$, then the claim holds trivially. We may thus assume $\| x\|>0$, so that we can set $\tilde{x}:=x/\| x\|$ and $\tilde{y}:=
y/\| x\|$ to get $\|\tilde{x}\|=1$ and~$\|\tilde{y}\|\le 1$. Due to~$\lambda\geq 1/2$, the claim holds if we have $\|x-y\|\leq\|x\|$. In the remaining case we have $\|\tilde x-\tilde y\|>1$, so that uniform convexity yields $\|\tilde x+\tilde y\|/2\leq 1-\eta(1)=\lambda$ and hence
\[ \|x+y\|=\|\tilde x+\tilde y\|\cdot\|x\|\leq 2\lambda\cdot\|x\|,\]
as required.\qed\\[1ex]
In order to connect with an equivalent definition of Rademacher type, we will use the following:

\begin{proposition}[Kahane-Kintchine inequality] 
For $q>1$ there exists a constant $K_q$ such that all finite 
sequences $(x_1,\ldots,x_n)$ in any Banach space $X$ validate
\[ \int^1_0 \left\|\sum^n_{i=1}r_i(t)x_i\right\| dt\le 
\left(\int^1_0 \left\|\sum^n_{i=1}r_i(t)x_i\right\|^q dt\right)^{1/q}\le 
K_q\cdot\int^1_0 \left\|\sum^n_{i=1}r_i(t)x_i\right\| dt. \]
Indeed, we may take $K_q:=\big((2q-1)/(q-1)\big)^{q-1}.$
\end{proposition}
{\bf Proof:} The result coincides with Theorem~1.e.13 of~\cite{Lindenstrauss-Tzafriri} (where $K_q$ is given in the proof).\hfill $\Box$ 
\\[1mm]
Pisier in~\cite{pisier73} uses a different definition of Rademacher type~$q$, which demands a constant~$c_q$ with
\[ \tag{$\ast\ast$} \ \EE\left(\left\| \sum^{n}_{i=1} \eps_ix_i\right\|^2\right)^{1/2} \le c_q\cdot\left(\sum^n_{i=1} 
\| x_i\|^q\right)^{1/q}. \]
Using the Kahane-Kintchine inequality, one sees that $(\ast)$ for $C_q$ entails $(\ast\ast)$ for $c_q:=K_2\cdot C_q$. In the converse direction, one can keep the constant and does not need Kahane-Kintchine:
\begin{lemma}\label{constants}
\begin{enumerate}
\item 
If $(\ast)$ holds for $C_q,$ then $(\ast\ast)$ holds for $c_q:=K_2C_q.$
\item 
If $(\ast\ast)$ holds for $c_q$, then $(\ast)$ holds for $C_q:=c_q$, where~$q\in [1,2]$ can be arbitrary.
\end{enumerate}
\end{lemma}
{\bf Proof:} 
1) Given $(\ast),$ we get from the Kahane-Kintchine inequality
\begin{multline*}
\mathbb E\left(\left\|\sum_{i=1}^n\varepsilon_ix_i\right\|^2\right)^{1/2}=\left(\int^1_0\left\|\sum_{i=1}^n r_i(t)x_i\right\|^2 dt\right)^{1/2}\leq K_2\cdot\int^1_0\left\|\sum^n_{i=1}r_i(t)x_i\right\| dt\leq\\ \leq 
K_2\cdot \left(\int_0^1\left\|\sum^n_{i=1} r_i(t)x_i\right\|^q dt\right)^{1/q}
=K_2\cdot \EE\left(\left\|\sum^n_{i=1}\eps_ix_i\right\|^q\right)^{1/q}\leq K_2\cdot C_q\cdot\left(\sum_{i=1}^n\|x_i\|^q\right)^{1/q}.
\end{multline*}
2) 
Given~$(\ast\ast)$, we get
\begin{multline*}
\EE\left(\left\|\sum^n_{i=1}\eps_ix_i\right\|^q\right) = \int^1_0\left\| \sum^n_{i=1}
r_i(t)x_i\right\|^qdt \le \left( \int^1_0\left\| \sum^n_{i=1}
r_i(t)x_i\right\|^2dt\right)^{q/2}
={}\\
{}=\left(\EE\left\|\sum^n_{i=1}\eps_ix_i\right\|^2\right)^{q/2}
\le c^q_q \cdot \sum^n_{i=1}\| x_i\|^q,
\end{multline*}
as $\left(\int^1_0 |f(t)|^q dt\right)^{1/q}\le
\left(\int^1_0 |f(t)|^2 dt\right)^{1/2}$ holds for square-integrable~$f$ and $q\in[1,2]$.\qed
\\[1ex]
By analysing the proof of Corollary 1 in \cite{pisier73}, we obtain 
the following numerical estimate:
\begin{theorem}\label{Rademacher}
Let $(X,\|\cdot\|)$ be a uniformly nonsquare Banach space with $\delta 
\in (0,1)$ witnessing this property. Define $\lambda:=1-\delta.$
Assume that $\xi\in(0,1)$ is so small and that $p'\in [2,\infty)$ is so large that
\[ \frac{1-\xi}{1+2\sqrt{2\xi}}\geq\frac{1}{2}\sqrt{2\lambda^2+2}\qquad\text{and}\qquad \frac{1}{2^{1/p'}}\geq 1-\xi.\]
Take $p$ with $1=\frac1{p}+\frac1{p'}$. Then for any $q\in(1,p)$, the space $X$ has Rade\-macher type~$q$ with constant 
\[ C_q=3\cdot\frac{2^{1/q}}{2^{(1/q)-(1/p)}-1}. \] 
\end{theorem}
{\bf Proof:} 
For $n\in\mathbb N$, Pisier~\cite{pisier73} defines $\mu_n(X)$ as the least real $\mu\geq 0$ such that
\[ \left(\int^1_0 \left\| \sum^n_{i=1} r_i(t)x_i\right\|^2dt
\right)^{1/2} \le \mu \cdot n\cdot\max_{1\le i\le n}
\| x_i\| \] 
holds for any finite sequence $(x_1,\ldots,x_n)$ in $X$. Similarly, he defines $\nu_n(X)$ as the least $\nu\geq 0$ with
\[ \left(\int^1_0 \left\| \sum^n_{i=1} r_i(t)x_i\right\|^2dt
\right)^{1/2} \le \nu \cdot \sqrt{n}\cdot
\left( \sum^n_{i=1}\| x_i\|^2\right)^{1/2}. \] 
The $\lambda$ of the present theorem is an upper bound for $\lambda_2(X)$ from~\cite{pisier73}, as noted after Lemma~\ref{lem:convex-to-nonsquare} above. Together with inequality~(2) on page~VII.10 of~\cite{pisier73} (corrected with the missing factor $1/n$), we get
 \[ \mu_2(X)\leq\frac{1}{2}\left[ \frac{4\lambda^2+4}{2}\right]^{1/2} 
=\frac{1}{2}\sqrt{2\lambda^2+2}\in (0,1). \]
For $\xi$ and~$p'$ as in the present theorem, the argument on pages~VII.10-11 of~\cite{pisier73} does now yield
\begin{equation*}
\nu_2(X)\leq 1-\xi\leq\frac1{2^{1/p'}}.
\end{equation*}
By Lemma~4 of~\cite{pisier73} it follows that $X$ has Rademacher type $q$ with a suitable constant~$C_q$, for any~$q$ as in the present theorem. To determine~$C_q$, we work out the proof of the cited lemma (which Pisier describes as analogous to the one of Lemma~2 from~\cite{pisier73}). Assume we have $\nu_N(X)\leq 1/N^{1/p'}$ for some integer $N\geq 2$ and real $p'\geq 2$ (we only need $N=2$ but state the original more general result). Given $q\in(1,p)$ with $1=(1/p)+(1/p')$, we shall establish $(\ast\ast)$ for suitable~$c_q$. Due to Lemma~\ref{constants}, it will follow that the Rademacher property~$(\ast)$ holds for~$C_q:=c_q$. Aiming at~$(\ast\ast)$, we consider an arbitrary sequence $(x_1,\ldots,x_n)$ in~$X$. For~$k\in\mathbb N$ we put
\[ A(k):=\left\{ j\in \{ 1,\ldots,n\}:\left( 
\frac{\sum^n_{i=1} \| x_i\|^q}{N^{k+1}}\right)^{1/q} <\| x_j\| 
\le  \left( 
\frac{\sum^n_{i=1} \| x_i\|^q}{N^{k}}\right)^{1/q}\right\}.\]
Write $|A(k)|$ for the cardinality of~$A(k)$. We pick a bijection~$f:\{1,\ldots,|A(k)|\}\to A(k)$ and compute
\begin{equation*}
\left(\int_0^1\left\|\sum_{i\in A(k)} r_i(t)x_i\right\|^2 dt\right)^{1/2}=\left(\int_0^1\left\|\sum_{i=1}^{|A(k)|} r_i(t)x_{f(i)}\right\|^2 dt\right)^{1/2}\leq\mu_{|A(k)|}(X)\cdot|A(k)|\cdot\max_{i\in A(k)}\|x_i\|.
\end{equation*}
Here the inequality holds by the definition of $\mu_n(X)$. The equality relies on the fact that the Rade\-macher functions represent independent copies of the same random variable, which allows us to omit the index shift from~$r_i$ to~$r_{f(i)}$. For each $i\in\{1,\ldots,n\}$ with $x_i\neq 0$, we have~$i\in A(k)$ for a unique integer~$k\geq 0$ (so that almost all $A(k)$ are empty). Using the Kahane-Kintchine inequality, we get
\begin{multline*}
\mathbb E\left(\left\|\sum_{i=1}^n\varepsilon_ix_i\right\|^2\right)^{1/2}=\left(\int^1_0\left\|\sum_{i=1}^n r_i(t)x_i\right\|^2 dt\right)^{1/2}\leq K_2\cdot\int^1_0\sum_{k=0}^{\infty}\left\|\sum_{i\in A(k)}r_i(t)x_i\right\| dt\leq\\
K_2\cdot\sum_{k=0}^\infty\left(\int_0^1\left\|\sum_{i\in A(k)} r_i(t)x_i\right\|^2 dt\right)^{1/2}\leq K_2\cdot\left(\sum_{k=0}^\infty \frac{\mu_{|A(k)|}(X)\cdot|A(k)|}{N^{k/q}}\right)\cdot\left(\sum_{i=1}^n\|x_i\|^q\right)^{1/q},
\end{multline*}
which is already close to~$(\ast\ast)$. We have $|A(k)|\leq N^{k+1}$, as on page~VII.5 of~\cite{pisier73}. By Proposition~3 and Lemma~3 from the same reference, this yields the first inequality in
\begin{equation*}
\mu_{|A(k)|}(X)\cdot|A(k)|\leq\nu_N(X)^{k+1}\cdot N^{k+1}\leq\frac{N^{k+1}}{N^{(k+1)/p'}}=N^{(k+1)/p}.
\end{equation*}
The second inequality and equality rely on the assumptions $\nu_N(X)\leq 1/N^{1/p'}$ and $1=(1/p)+(1/p')$. One can conclude
\begin{equation*}
K_2\cdot\sum_{k=0}^\infty \frac{\mu_{|A(k)|}(X)\cdot|A(k)|}{N^{k/q}}\leq K_2\cdot\sum_{k=0}^\infty\frac{N^{1/p}}{N^{k\cdot(1/q-1/p)}}=K_2\cdot\frac{N^{1/q}}{N^{(1/q)-(1/p)}-1}=:c_{N,q}.
\end{equation*}
We have thus established $(\ast\ast)$ with $c_{N,q}$ at the place of~$c_q$. Under the assumptions of the theorem, we get the Rademacher property~$(\ast)$ for $C_q:=c_{2,q}$, as seen above. Now we only need to note that $K_2=3$.\qed
\\[1mm]
In the next section we will need the following consequence 
of the Rademacher property:
\begin{proposition}[{\cite[Proposition 9.11]{ledoux-talagrand}}]\label{condition}
Assume that $(X,\|\cdot\|)$ is a Banach space of Rademacher type $q\in [1,2]$
with constant $C_q.$ Then, for every finite sequence $X_1,\ldots,X_n$ of 
independent mean zero Radon random variables in 
$L_q(X)$, one has
\[\tag{$+$} \EE\left(\left\|\sum_{i=1}^n X_i\right\|^q\right)^{1/q}\le 2C_q\cdot\left(\sum_{i=1}^n \EE\left(\| X_i\|^q\right)\right)^{1/q}.\]
\end{proposition}
Together with the previous results, the proposition yields $c$ and $q$ as in condition~(9) of~\cite{FreundKohlenbach}.
\begin{remark} \rm 
We will only need the special case of Proposition~\ref{condition} in which each of the $X_i$ assumes value $(x_j-x)/q$ with probability~$\lambda_j$, for a given convex combination $x=\sum\lambda_jx_j$. In this case, one can deduce the proposition by elementary manipulations of finite sums.
\end{remark}

\section{A rate of asymptotic regularity for ergodic averages in 
uniformly convex Banach spaces}\label{section-Bruck}
For Hilbert spaces, there is an easy quadratic rate of asymptotic regularity 
for the sequence of Ces\`aro means, which is given (see e.\,g.~\cite{BrezisBrowder} and~\cite{Kohlenbach(baillon)}) by
\[ \left\| \frac{1}{n}\sum^{n-1}_{i=0} T^ix-
T\left(  \frac{1}{n}\sum^{n-1}_{i=0} T^ix\right) \right\|\le 
\frac{1}{\sqrt{n}}\cdot\operatorname{diam}(C)\quad\text{for}\quad n\in\mathbb N\backslash\{0\}.\]
The asymptotic regularity for the Ces\`aro means in uniformly convex 
Banach spaces was first proved in \cite{Bruck(81)}. Although this proof by Bruck is essentially constructive, the concrete rate of convergence 
hidden in the proof is left implicit. In this section we extract the explicit rate that is specified in the following theorem. Let us emphasize that $q$ and~$C_q$ can be chosen according to Theorem~\ref{Rademacher}, where we may take $\delta=\eta(1)$ due to Lemma~\ref{lem:convex-to-nonsquare}. Hence our rate depends on~$\eta$ and~$b$ only, where $\eta$ is a modulus of uniform convexity of the 
space in question and $b>0$ is such that $C\subseteq B_{b/2}(0):=
\{ x\in X: \| x\|\le b/2\}.$ We point out that the superscript of~$\xi$ denotes iterations, which are explained by $\xi^0(t):=t$ and $\xi^{p+1}(t):=\xi(\xi^p(t))$.

\begin{theorem}\label{thm:asymptotic-reg}
Let $(X,\|\cdot\|)$ be a Banach space that is uniformly convex with modulus $\eta$. Consider a nonexpansive map $T:C\to C$ on a nonempty subset $C\subseteq B_{b/2}(0)$ that is closed and convex. Assume that $X$ is of Rademacher type $q\in (1,2]$ with constant $C_q$. Given $\varepsilon>0$, pick $\tilde p\in\mathbb N$ so large that we have $2C_q\cdot\tilde p^{(1-q)/q}\leq\varepsilon/(9b)$. Consider $p\in\mathbb N$ with $p\geq 2b/\delta^2$ for
\begin{equation*}
\delta:=\xi^{\tilde p}\left(\frac{\varepsilon}{9}\right)\quad\text{with}\quad\xi(t):=\frac{t}{12}\cdot\eta\left(\min\left\{2,\frac{t}{b}\right\}\right).
\end{equation*}
For any~$\alpha<\xi^{p-1}(\delta^2/2)$ with $0<\alpha<\varepsilon/3$ and arbitrary $x\in C$, we then have
\[  \left\| \frac{1}{n}\sum^{n-1}_{i=0} T^ix-
T\left(  \frac{1}{n}\sum^{n-1}_{i=0} T^ix\right)\right\|\leq\eps\quad\text{for all}\quad n\ge\frac{b}{\alpha}.\]
\end{theorem}

The proof of the theorem will occupy us for the rest of this section. As noted before, we use a modulus $\eta:(0,2]\to(0,1]$ that satisfies the condition in Lemma~\ref{lem:convex-to-nonsquare} but need not be optimal. We set
\begin{equation*}
\eta_1(0):=0\quad\text{and}\quad\eta_1(\varepsilon):=\sup\left\{\eta(\varepsilon')\,|\,0<\varepsilon'\leq\min\{2,\varepsilon\}\right\}
\end{equation*}
to get a function~$\eta_1:[0,\infty)\to[0,1]$ that is increasing. Now let $\tilde\eta:[0,\infty)\to[0,\infty)$ be given by
\begin{equation*}
\tilde{\eta}(\eps):=\frac{1}{2}\cdot\int^{\eps}_0\eta_1(t)dt.
\end{equation*}
The point is that this makes~$\tilde\eta$ convex. Given that $\varepsilon\in(0,2]$ entails $0<\tilde\eta(\varepsilon)\leq\eta_1(\varepsilon)$, we see that $\tilde\eta$ is still a modulus of uniform convexity for~$X$. We now consider the function $\gamma:[0,\infty)\to [0,\infty)$ with
\[ \gamma(\eps):=\frac{b}{2}\cdot\tilde{\eta}\left(\frac{4\eps}{b}\right).\]
Let us note that $\gamma$ is continuous and strictly increasing with image~$[0,\infty)$. We thus get a continuous and strictly increasing inverse~$\gamma^{-1}:[0,\infty)\to[0,\infty)$. Furthermore, the function~$\gamma$ is convex. Our functions $\tilde\eta$ and $\gamma$ coincide with~$\eta$ and $\gamma_2$ from Remark~2.3 and Definition~2.1 of~\cite{FreundKohlenbach}, respectively (given that $\eta_1(t)\leq 1$ entails $\tilde\eta(\eps)\leq\eps/2$ and hence $\gamma(\varepsilon)\leq\varepsilon$). By~\cite[Lemma~2.2]{FreundKohlenbach} (essentially a result of Bruck~\cite{Bruck79}), we learn that $T$ is of type~($\gamma$), i.\,e., that 
\[ \gamma(\| T(\lambda x_1+(1-\lambda)x_2)-(\lambda Tx_1+(1-\lambda)Tx_2)\|)
\le \| x_1-x_2\|-\| Tx_1-Tx_2\|\]
holds for any $x_1,x_2\in C$ and all $\lambda\in [0,1]$.

\begin{lemma}\label{lem:tilde-q}
For $\tilde q:[0,\infty)\to [0,\infty)$ with $\tilde{q}(\eps):=\gamma^{-1}(3\eps)+\eps$ and for~$\xi$ as in Theorem~\ref{thm:asymptotic-reg}, we have
\begin{equation*}
t\in\big[0,\xi^p(\varepsilon)\big)\quad\Rightarrow\quad\tilde q^p(t)\in[0,\varepsilon).
\end{equation*}
\end{lemma}
{\bf Proof:}
We use induction to reduce to the case of~$p=1$. The claim is immediate for~$p=0$. In the induction step, the case of~$p=1$ (with $\xi^p(\varepsilon)$ at the place of~$\varepsilon$) and the induction hypothesis yield
\begin{equation*}
t\in\big[0,\xi^{p+1}(\varepsilon)\big)=\big[0,\xi(\xi^p(\varepsilon))\big)\quad\Rightarrow\quad\tilde q(t)\in\big[0,\xi^p(\varepsilon)\big)\quad\Rightarrow\quad\tilde q^{p+1}(t)=\tilde q^p(\tilde q(p))\in[0,\varepsilon).
\end{equation*}
To establish the result for~$p=1$, we first observe that we have
\[ \tilde{\eta}(\eps)\ge
\frac{1}{2}\cdot\frac{\eps}{2}\cdot\eta_1\left(\frac{\eps}{2}\right)\ge
\frac{\eps}{4}\cdot\eta\left(\min\left\{2,\frac{\eps}{2}\right\}\right),\] 
as $\eta_1$ is increasing. With $2\varepsilon/b$ at the place of~$\varepsilon$, we get
\begin{equation*}
\xi(\varepsilon)=\frac{\varepsilon}{12}\cdot\eta\left(\min\left\{2,\frac{\varepsilon}{b}\right\}\right)\leq\frac{b}{6}\cdot\tilde\eta\left(\frac{2\varepsilon}{b}\right)=\frac13\cdot\gamma\left(\frac{\varepsilon}{2}\right).
\end{equation*}
Since $\tilde q$ is strictly increasing, it follows that $t<\xi(\varepsilon)$ entails
\begin{equation*}
\tilde q(t)<\tilde q\left(\frac13\cdot\gamma\left(\frac{\varepsilon}{2}\right)\right)=\gamma^{-1}\left(\gamma\left(\frac{\varepsilon}{2}\right)\right)+\frac13\cdot\gamma\left(\frac{\varepsilon}{2}\right)<\varepsilon.
\end{equation*}
Here the last inequality holds because we have $\gamma(t)\leq t$, as observed above.\hfill\qed\\[1mm]
The following is a final preparation for the proof of our main theorem.
\begin{lemma}\label{lem:q_n}
For $n\in\mathbb N\backslash\{0\}$, we consider $q_n:[0,\infty)\to[0,\infty)$ with $q_n(\varepsilon):=\gamma^{-1}(2\varepsilon+(b/n))+\varepsilon$. Given~$\varepsilon>0$, we get
\begin{equation*}
n\geq\frac{b}{\varepsilon}\quad\Rightarrow\quad q_n^p(\varepsilon)\leq\tilde q^p(\varepsilon)\text{ for all $p\in\mathbb N$ }.
\end{equation*}
\end{lemma}
{\bf Proof:}
Since $n\geq b/\varepsilon$ entails $2\varepsilon+(b/n)\leq 3\varepsilon$, the result for~$p=1$ follows from the fact that $\gamma^{-1}$ is increasing. We derive the general case by induction on~$p$. As before, the case of~$p=0$ is immediate. In the induction step, the induction hypothesis and the fact that $q_n$ is increasing yield
\begin{equation*}
q_n^{p+1}(\varepsilon)=q_n(q_n^p(\varepsilon))\leq q_n(\tilde q^p(\varepsilon)).
\end{equation*}
Since $\tilde q$ is increasing with $\varepsilon\leq\tilde q(\varepsilon)$, we get $\varepsilon\leq\tilde q^p(\varepsilon)$ by an auxiliary induction on~$p$. Given $n\geq b/\varepsilon$, we thus have $n\geq b/\tilde q^p(\varepsilon)$. Now the result for~$p=1$ (with $\tilde q^p(\varepsilon)$ at the place of~$\varepsilon$) yields
\begin{equation*}
q_n(\tilde q^p(\varepsilon))\leq \tilde q(\tilde q^p(\varepsilon))=\tilde q^{p+1}(\varepsilon),
\end{equation*}
as needed to complete the induction step.\hfill\qed\\[1mm]
We now have all ingredients to show the result that was stated above:\\[1mm]
{\bf Proof of Theorem~\ref{thm:asymptotic-reg}:}
Let $F_\delta(T)=\{x\in C\,:\,\|x-Tx\|\leq\delta\}$ be the set of $\delta$-approximate fixed points of~$T$. For $S\subseteq X$, we write
\begin{equation*}
\co_p(S):=\left\{\left.\sum_{i=1}^n\lambda_ix_i\,\right|\,x_i\in S\text{ and }\lambda_i\geq 0\text{ with }\sum_{i=0}^n\lambda_i=1\text{ for }n\leq p\right\}
\end{equation*}
for the set of convex combinations of at most $p$ elements. By $\co(S):=\bigcup_{p\in\mathbb N}\co_p(S)$ we denote the convex hull. For arbitrary~$\varepsilon>0$, we will show that the $\delta$ that is specified in Theorem~\ref{thm:asymptotic-reg} validates
\begin{equation}\label{eq:clco-F_delta}\tag{$++$}
\co(F_{\delta}(T))\subseteq F_{\eps/3}(T).
\end{equation}
Before we prove this, we show how to deduce the theorem. First note that the closure of $\co(F_{\delta}(T))$ will still be contained in $F_{\eps/3}(T)$, as the latter is closed. As in Theorem~\ref{thm:asymptotic-reg}, we assume $p\geq 2b/\delta^2$ and~$\alpha<\xi^{p-1}(\delta^2/2)$. For $\tilde q$ and $q_n$ as above, Lemmas~\ref{lem:tilde-q} and~\ref{lem:q_n} yield
\begin{equation*}
q_n^{p-1}(\alpha)\leq \tilde q^{p-1}(\alpha)<\frac{\delta^2}{2}\quad\text{for all}\quad n\geq\frac{b}{\alpha}.
\end{equation*}
For $x_n:=T^n x$ we clearly get $\|x_{n+1}-Tx_n\|=0\leq\alpha$. We have established all properties that are used in the proof of Theorem~1.3 from~\cite{Bruck(81)}, which shows that any $x\in C$ validates
\begin{equation*}
\frac1n\cdot\sum_{i=0}^{n-1}T^ix\in F_\varepsilon(T)\quad\text{for all}\quad n\geq\frac{b}{\alpha}.
\end{equation*}
This coincides with the conclusion of the desired Theorem~\ref{thm:asymptotic-reg}. It remains to show that the~$\delta$ from the theorem satisfies~(\ref{eq:clco-F_delta}). Due to the assumption that $X$ has Rademacher type~$q$ with constant~$C_q$, we can apply Proposition~\ref{condition}. As noted above, this yields condition~(9) of~\cite{FreundKohlenbach}, with $X$ and $2C_q$ at the place of~$X^2$ and $c$. Essentially by~\cite[Lemma~2.6]{FreundKohlenbach} (based on the proof of~\cite[Theorem~1.1]{Bruck(81)}), we get
\begin{equation*}
\co(M)\subseteq\co_{\tilde p}(M)+B_{\varepsilon/9}\quad\text{for any}\quad M\subseteq B_{b/2},
\end{equation*}
provided we have $2C_q\cdot\tilde p^{(1-q)/q}\leq\varepsilon/(9b)$ as in Theorem~\ref{thm:asymptotic-reg}. The function~$\tilde q$ is continuous and strictly increasing with image~$[0,\infty)$, as the same holds for~$\gamma$. Let $\sigma=\tilde q^{-1}:[0,\infty)\to[0,\infty)$ be its inverse. As in the proof of Theorem 1.2 in \cite{Bruck(81)}, one can show 
\[\delta\le\sigma^{\tilde{p}}(\eps/9)\quad\Rightarrow\quad 
\co(F_{\delta}(T))\subseteq F_{\eps/3}(T).\]
To obtain~(\ref{eq:clco-F_delta}) for~$\delta=\xi^{\tilde p}(\varepsilon/9)$ as in Theorem~\ref{thm:asymptotic-reg}, we show $\xi^i(t)\leq\sigma^i(t)$ by induction on~$i\in\mathbb N$. Let us first note that Lemma~\ref{lem:tilde-q} yields~$\tilde q(\xi(s))\leq s$, as $\tilde q$ is continuous. Since $\sigma=\tilde q^{-1}$ is increasing, we can conclude $\xi(s)\leq\sigma(s)$. Given $s:=\xi^i(t)\leq\sigma^i(t)$, we thus get
\begin{equation*}
\xi^{i+1}(t)=\xi(\xi^i(t))\leq\sigma(\xi^i(t))\leq\sigma(\sigma^i(t))=\sigma^{i+1}(t),
\end{equation*}
as needed for the induction step.\hfill\qed\\[1mm]
As promised in the introduction, we now discuss an alternative rate of asymptotic regularity:
\begin{remark}\rm
As shown by Bruck (see~\cite[Theorem~2.1]{Bruck(81)}), the above function~$\gamma$ can be transformed into a continuous, strictly increasing and convex $\tilde\gamma:[0,\infty)\to[0,\infty)$ with $\tilde\gamma(0)=0$ such that
\begin{equation*}
\tilde\gamma\left(\left\|T\left(\sum_{i=1}^n\lambda_ix_i\right)-\sum_{i=1}^n\lambda_i Tx_i\right\|\right)\leq\max_{1\leq i,j\leq n}\left(\|x_i-x_j\|-\|Tx_i-Tx_j\|\right)
\end{equation*}
holds for any convex combination $\sum\lambda_ix_i$. By a result of Zhu, Huang and Li, the conclusion of our Theorem~\ref{thm:asymptotic-reg} holds for all~$n$ above a certain $n_0$ that depends on such a~$\tilde\gamma$ (see Lemma~3.5 of~\cite{zhu}, which works in a more general semigroup setting). In Section~2 of~\cite{FreundKohlenbach}, we have shown how $\tilde\gamma$ can be expressed in terms of $b,c,q$ and $\eta$, assuming that $X^2$ has Rademacher type~$q$ with constant~$c$ (but note that the $\tilde\gamma$ in~\cite{FreundKohlenbach} is not made convex). Now Section~\ref{section-Pisier} of the present paper shows how to express $c$ and $q$ in terms of~$\eta$. This makes it possible to express the rate of Zhu, Huang and Li in terms of~$\eta$ and~$b$. A very rough comparison suggests that our rate from Theorem~\ref{thm:asymptotic-reg} is better, as it involves fewer iterations of~$\eta$ when~$q$ is close to~$1$. However, we have not established a precise comparison between the two rates.
\end{remark}

To conclude, we observe how certain assumptions can be weakened:

\begin{remark}\label{rmk:zhu}\rm
First, note that the proof above involves $x_n=T^nx$ with $\| x_{n+1}-Tx_n\|=0$. The argument does also go through for different $x_n$ that satisfy $\| x_{n+1}-Tx_n\|\le\alpha$ for all~$n$. Secondly, the assumption that $C\subseteq B_{b/2}(0)$ is bounded can be secured when $T:C'\to C'$ is defined on an unbounded set $C'$ (still assumed to be closed and convex) and has a fixed point~$f=T(f)$, as pointed out in~\cite{FreundKohlenbach}. Indeed, for given $x\in C'$ we can then consider
\begin{equation*}
C:=C'\cap B_{\|x-f\|}(f)\subseteq B_{\|x-f\|+\|f\|}(0),
\end{equation*}
which is closed and convex with $T(C)\subseteq C\ni x$. In fact, it is enough to have a $d>0$ such that $T$ has arbitrarily good approximate fixed points in $B_d(0)$, i.\,e., such that $B_d(0)\cap F_\varepsilon(T)\neq\emptyset$ holds for all~$\varepsilon>0$. Indeed, for any $f\in F_\varepsilon(T)$ we inductively get
\begin{equation*}
\|T^nf-f\|\leq\|T^nf-T^{n-1}f\|+\|T^{n-1}f-f\|\leq\|Tf-f\|+(n-1)\cdot\varepsilon\leq n\cdot\varepsilon.
\end{equation*}
Given $n\in\mathbb N\backslash\{0\}$ and for arbitrary~$\varepsilon>0$, we now pick $f\in B_d(0)\cap F_{\varepsilon/n}(0)$ to obtain
\begin{equation*}
\|T^nx\|\leq\|T^nx-T^nf\|+\|T^nf-f\|+\|f\|\leq\|x-f\|+\varepsilon+\|f\|\leq\|x\|+2\cdot d+\varepsilon.
\end{equation*}
As~$\varepsilon>0$ was arbitrary, we may now omit it, for each~$n\in\mathbb N$. This allows us to consider
\begin{equation*}
R:=\limsup_{n\to\infty}\|x-T^nx\|\leq 2\cdot (\|x\|+d).
\end{equation*}
As in the proof of~\cite[Theorem~1]{Reich-bounded-orbit} (for the simple case where $(a_{nk})$ is the identity matrix), the set
\begin{equation*}
C:=\left\{y\in C'\,\left|\,\limsup_{n\to\infty}\|y-T^nx\|\leq R\right.\right\}\subseteq B_{R+\|x\|+2d}(0)
\end{equation*}
is closed and convex with $T(C)\subseteq C\ni x$. Note that the Browder-G\"ohde-Kirk fixed point theorem will now yield an actual fixed point in~$C$. The fact that approximate fixed points can play the role of actual ones is also guaranteed by logical metatheorems from proof mining (see Corollaries~5.2 and~6.8 as well as the comment following Remark~5.9 of~\cite{KohlenbachGerhardy}).
\end{remark}

\noindent
{\bf Acknowledgment:} Both authors were supported by the `Deutsche Forschungs\-gemein\-schaft' (DFG, German Research Foundation) -- Projects 460597863 and DFG KO 1737/6-2.

\end{document}